\documentclass[12pt]{amsart}
\usepackage{latexsym, amssymb, amsfonts, amscd}

\theoremstyle{plain}
\newtheorem{Thm}{Theorem}[subsection]
\newtheorem{Cor}[Thm]{Corollary}

\newtheorem{Prop}[Thm]{Proposition}
\newtheorem{Lem}[Thm]{Lemma}

\newtheorem{Thm'}{Theorem}[section]
\newtheorem{Cor'}[Thm']{Corollary}
\newtheorem{Prop'}[Thm']{Proposition}
\newtheorem{Lem'}[Thm']{Lemma}
\newtheorem{Cl'}[Thm']{Claim}

\theoremstyle{definition}
\newtheorem{Def}[Thm]{Definition}
\newtheorem{Rem}[Thm]{Remark}

\newtheorem{Emp}[Thm]{}

\newtheorem{Not}[Thm]{Notation}

\newtheorem{Def'}[Thm']{Definition}
\newtheorem{Rem'}[Thm']{Remark}
\newtheorem{Rem1'}[Thm']{Remarks}
\newtheorem{Emp'}[Thm']{}
\newtheorem{Ex'}[Thm']{Example}
\newtheorem{Exs'}[Thm']{Examples}
\newtheorem{Con'}[Thm']{Construction}
\newtheorem{Not'}[Thm']{Notation}
\newtheorem{Q'}[Thm']{Question}

\numberwithin{equation}{section}

\newcommand{\qlbar}{\overline{\B{Q}_l}}
\newcommand{\fqbar}{\overline{\fq}}

\newcommand{\La}{\Lambda}
\newcommand{\lan}{\left\langle}
\newcommand{\ran}{\right\rangle}
\newcommand{\pp}{\boxtimes}\newcommand{\ov}{\overline}
\newcommand{\un}{\underline}

\newcommand{\fq}{\B{F}_q}

\newcommand{\B}[1]{\mathbb#1}
\newcommand{\cal}[1]{\mathcal{#1}}
\newcommand{\form}[1]{(\ref{Eq:#1})}
\newcommand{\C}[1]{\cal#1}

\newcommand{\isom}{\overset {\thicksim}{\to}}

\newcommand{\lra}{\longrightarrow}
\newcommand{\lla}{\longleftarrow}
\newcommand{\hra}{\hookrightarrow}
\newcommand{\wt}{\widetilde}
\newcommand{\diez}{\natural}

\newcommand{\Gm}{\Gamma}

\newcommand{\Dt}{\Delta}

\newcommand{\rs}[1]{Section \ref{S:#1}}
\newcommand{\rl}[1]{Lemma \ref{L:#1}}

\newcommand{\rp}[1]{Proposition \ref{P:#1}}
\newcommand{\rr}[1]{Remark \ref{R:#1}}

\newcommand{\re}[1]{\ref{E:#1}}
\newcommand{\rco}[1]{Corollary \ref{C:#1}}

\newcommand{\rt}[1] {Theorem \ref{T:#1}}

\newcommand{\sm}{\smallsetminus}

\newcommand{\Spec}{\operatorname{Spec}}

\newcommand{\Tr}{\operatorname{Tr}}
\newcommand{\Fr}{\operatorname{Fr}}

\newcommand{\Id}{\operatorname{Id}}

\newcommand{\Hom}{\operatorname{Hom}}

\begin{document}

\title[A proof of a generalization of Deligne's conjecture]%
{A proof of a generalization of Deligne's conjecture}
\author{Yakov Varshavsky}
\address{Institute of Mathematics\\
Hebrew University\\
Givat-Ram, Jerusalem,  91904\\
Israel}
\email{vyakov@math.huji.ac.il }

\thanks{The work was supported by 
THE ISRAEL SCIENCE FOUNDATION (Grant No. 555/04)}
\date{August 2005}

\begin{abstract}

The goal of this paper is to give a simple proof of Deligne's conjecture on 
the Lefschetz trace formula (proven by Fujiwara) 
and to generalize it to the situation appearing in \cite{KV}. 
Our proof holds in the realm of ordinary algebraic geometry and does not use 
rigid geometry.
\end{abstract}
\maketitle

\section*{Introduction}

Suppose we are given a correspondence $X\overset{c_1}{\lla}C\overset{c_2}{\lra}X$
of separated schemes of finite type over a separably closed field $k$, 
an ``$\ell$-adic sheaf'' $\C{F}\in D^b_{c}(X,\qlbar)$ and a morphism $u:c_{2!}c_1^*\C{F}\to\C{F}$. 
If $c_1$ is proper, then $u$ gives rise to an endomorphism 
$R\Gm_c(u):R\Gm_c(X,\C{F})\to R\Gm_c(X,\C{F})$. 

When $X$ is proper, the general Lefschetz--Verdier trace formula
\cite[Cor. 4.7]{Il} asserts that the trace $\Tr(R\Gm_c(u))$ equals the sum  
$\sum_{\beta\in\pi_0(Fix(c))}LT_{\beta}(u)$, where 
$Fix(c):=\{y\in C\,|\,c_1(y)=c_2(y)\}$ is the scheme of fixed points of $c$,
and $LT_{\beta}(u)$ is a so called ``local term'' of $u$ at $\beta$.
This result has two drawbacks: it fails when $X$ is not proper, and the 
``local terms'' are very inexplicit. 

Deligne conjectured that the situation becomes better if  
$X\overset{c_1}{\lla}C\overset{c_2}{\lra}X$ is defined over a 
finite field $\fq$,  $c_2$ is quasi-finite, and we twist $c_1$ by a sufficiently 
high power of the geometric Frobenius morphism. 
More precisely, he conjectured that in this case 
the Lefschetz--Verdier trace formula holds also for a non-proper $X$'s, 
$Fix(c)$ is finite, and for each $y\in c_1^{-1}(x)\cap c_2^{-1}(x)$, 
the local term $LT_{y}(u)$  equals the trace of the 
endomorphism $u_y:\C{F}_x\to\C{F}_x$, induced by $u$.

The conjecture was first proven by  
Pink \cite{Pi} assuming the resolution of singularities, and then by 
Fujiwara \cite{Fu} unconditionally. 

Theorem of Fujiwara has a fundamental importance for Langlands' program. 
For example, it was crucially used by Flicker--Kazhdan, Harris--Taylor, and Lafforgue. 

In a joint project \cite{KV} with David Kazhdan on the global Langlands correspondence 
over function fields we needed a generalization of Deligne's conjecture. Namely, 
instead of assuming that $c_1$ is proper we assumed that 
there exists an open subset $U\subset X$ such that $c_1|_{c_1^{-1}(U)}$ is proper, 
$X\sm U$ is ``locally $c$-invariant'', and $\C{F}$ vanishes on $X\sm U$. 
In this case, $u$ still gives rise to an endomorphism 
$R\Gm_c(u)$, and the main result of the present work asserts that the conclusion of 
Deligne's conjecture holds in this case. 

The strategy of our proof is similar to that of \cite{Pi} and \cite{Fu}: first we 
reduce the problem to vanishing of local terms $LT_{\beta}$, then 
we make the correspondence ``contracting'' by twisting it with a 
sufficiently high power of Frobenius, and finally we show the vanishing of local terms 
for ``contracting'' correspondences. 

Our approach differs from that of Fujiwara in two respects. First of all, 
our notion of a ``contracting'' correspondence is much simpler both to define and to use. 
Namely, we use the most naive notion of an ``infinitesimally'' contracting correspondence, 
which has a simple geometric description in terms of a ``deformation to the normal cone''. 
As a result, our bound on the power of Frobenius is sharper and more explicit.

Secondly, to prove a generalization of Deligne's conjecture described above, 
we work ``locally''. More precisely, as in \cite{Pi}, to show the vanishing of 
``local terms'', we first show the vanishing of so called ``trace maps'', from which 
``local terms'' are obtained by integration. 

Finally we would like to emphasize that our result is ultimately an assertion 
about the geometry of correspondences and not about sheaves. Also our proof applies without 
any changes to more general situations like compactifiable algebraic spaces or 
Deligne--Mumford stacks.

\section*{Notation and conventions}

Throughout the paper, all schemes will be separated schemes of 
finite type over a separably closed field $k$. We fix a prime $l$, invertible in $k$,  
and a commutative ring with identity $\La$, which is either finite and 
is annihilated by some power of $l$, or a finite extension of $\B{Z}_l$ or $\B{Q}_l$. 

To each scheme $X$ as above, we associate a derived category $D_{ctf}^b(X,\La)$  
of complexes of $\La$-modules of finite tor-dimension with constructible cohomology.  
This category is known to be stable under the six operations
$f^*, f^!, f_*, f_!, \otimes$ and $\C{RHom}$. 

For each $X$ as above, we denote by $\pi_X:X\to\Spec k$ the structure morphism, by 
$K_X=\pi_X^!\La$ the dualizing complex of $X$, and 
by $\B{D}=\C{RHom}(\cdot,K_X)$ the Verdier duality functor. We also write 
$R\Gm_c(X,\cdot)$ instead of $\pi_{X!}$. For an embedding 
$f:Y\hra X$, we will write $\cdot|_Y$ instead of $f^*$.

 For a closed subscheme $Z\subset X$, denote by $\C{I}_Z\subset\C{O}_X$ the sheaf of 
ideals of $Z$. For a morphism $f:X'\to X$, denote by $f^{-1}(Z)$ the schematic inverse 
image of $Z$, i.e., the closed subscheme of $X'$ such that 
$\C{I}_{f^{-1}(Z)}=f^{\cdot}(\C{I}_{Z})\cdot\C{O}_{X'}\subset
\C{O}_{X'}$, where $f^{\cdot}$ is the pullback map for functions.
We will also identify a closed subset of $X$ with the corresponding 
closed reduced subscheme.

Let $\fqbar$ be an algebraic closure of the finite field $\fq$. 
We say that an object $\C{X}$ over $\fqbar$ is {\em defined over} $\fq$, if it is 
a pullback of the corresponding object over $\fq$. 

\section*{Acknowledgments}
 
This work would not be possible without David Kazhdan, who among other things  
explained to me how to define $R\Gm_c(u)$ in the case described above and  
suggested that an analog of Deligne's conjecture should work in this situation.

Also the author thanks Alexander Beilinson, who gave a reference to \cite{Ve} and whose
comments and suggestions helped to simplify and improve the exposition. 

Part of the work was done while the author visited 
the University of Chicago and Northwestern University. The author thanks these institutions for 
hospitality and financial support.
\section{Formulation of the result} \label{S:def}

\subsection{Preliminaries on correspondences.}

\begin*
\vskip 8truept
\end*

\begin{Def} \label{D:corr} 
By a {\em correspondence}, we mean a morphism of schemes of the form
$c=(c_1,c_2):C\to X_1\times X_2$. By a {\em $c$-morphism}, we mean a morphism 
$u:c_{2!}c_1^*\C{F}_1\to\C{F}_2$ for some $\C{F}_1\in D^b_{ctf}(X_1,\La)$ and 
$\C{F}_2\in D^b_{ctf}(X_2,\La)$. 
\end{Def}

\begin{Rem}
A $c$-morphism is usually called a {\em cohomological correspondence lifting $c$}.
\end{Rem}

\begin{Emp} \label{E:restr}
{\bf Restriction of correspondences.}
 Let $c$ and $u$ be as in \ref{D:corr}, and 
$C^0\subset C$, $X^0_1\subset X_1$ and $X^0_2\subset X_2$ open subschemes 
(resp. closed subschemes satisfying $c_2^{-1}(X^0_2)_{red}=C^0_{red}$) such that 
$c$ induces a correspondence $c^0:C^0\to  X^0_1\times X^0_2$. 

Then we have a base change morphism 
$BC:c^0_{2!}(\C{F}|_{C^0})\to (c_{2!}\C{F})|_{X^0_2}$ for every $\C{F}\in  D^b_{ctf}(C,\La)$, 
hence $u$ restricts to a $c^0$-morphism 
\[
u^0:c^0_{2!}c^{0*}_1(\C{F}_1|_{X^0_1})=c^0_{2!}(c_1^*\C{F}_1|_{C^0})\overset{BC}{\lra}
(c_{2!}c_1^*\C{F}_1)|_{X^0_2}\overset{u}{\lra} \C{F}_2|_{X^0_2}.
\]
\end{Emp}

\begin{Emp} \label{E:rgm}
{\bf Endomorphism on the cohomology.}
Let $c$ and $u$ be as in \ref{D:corr}.

a) Assume that $c_1$ is proper. Then $u$ induces a morphism
\[
R\Gm_c(u):R\Gm_c(X_1,\C{F}_1)\overset{c_1^*}{\lra}R\Gm_c(C, c_1^*\C{F}_1)=
R\Gm_c(X_2,c_{2!}c_1^*\C{F}_1)\overset{u}{\lra}R\Gm_c(X_2,\C{F}_2).
\]

b) More generally, assume that there exists an open subset $U_1\subset X_1$ such that 
$c_1|_{c_1^{-1}(U_1)}:c_1^{-1}(U_1)\to U_1$ is proper and $\C{F}_1|_{X_1\sm U_1}=0$.
Let $c^0:c_1^{-1}(U_1)\to U_1\times X_2$ be the restriction of $c$, and 
$u^0$ the restriction of $u$ to $c^0$ (see \ref{E:restr}). Then $c^0_1$ is proper, 
hence by a) $u^0$  gives rise to a morphism 
$R\Gm_c(u^0): R\Gm_c(U_1,\C{F}_1|_{U_1})\to R\Gm_c(X_2,\C{F}_2)$.

On the other hand, the canonical map 
$j_!:R\Gm_c(U_1,\C{F}_1|_{U_1})\to R\Gm_c(X_1,\C{F}_1)$ is an isomorphism, therefore 
$u$ gives rise to a morphism
\[
R\Gm_c(u):=R\Gm_c(u^0)\circ(j_!)^{-1}:R\Gm_c(X_1,\C{F}_1)\to R\Gm_c(X_2,\C{F}_2).
\]
\end{Emp}

\begin{Emp}
{\bf Twisting of correspondences.}
For a scheme $X/\fqbar$, which defined over $\fq$, we denote by $\Fr_X=\Fr_{X,q}:X\to X$
the geometric Frobenius morphism. For a correspondence 
$c:C\to X_1\times X_2$ defined over $\fq$ and an integer $n\in\B{N}$, we denote by $c^{(n)}$
the correspondence $(c_1^{(n)},c_2):C\to X_1\times X_2$, 
where $c_1^{(n)}=\Fr_{X_1}^n\circ c_1=c_1\circ \Fr_C^n$.
\end{Emp}

\begin{Not} \label{N:partcase}
For a correspondence $c:C\to X\times X$, let $\Dt$ be the diagonal map 
$X\hra X\times X$, put $Fix(c):=c^{-1}(\Dt(X))\subset C$, and denote 
by $\Dt': Fix(c)\hra C$ and $c':Fix(c)\to X$ 
the inclusion map and the restriction of $c$, respectively.
We call $Fix(c)$ {\em the scheme of fixed points of $c$}.
\end{Not}

\begin{Emp} \label{E:naive}
{\bf Quasi-finite case.}
Let $c:C\to X\times X$ be a correspondence such that  $c_2$ is quasi-finite, and let  
$u:c_{2!}c_1^*\C{F}\to\C{F}$ be a $c$-morphism.

Each fiber $(c_{2!}c_1^*\C{F})_x$ decomposes as  
$\oplus_{y\in c_2^{-1}(x)}(c_1^*\C{F})_y=\oplus_{y\in c_2^{-1}(x)}\C{F}_{c_1(y)}$.
Hence for each $y\in c_2^{-1}(x)$, the map $u$ induces a morphism
$u_y:=u_x|_{\C{F}_{c_1(y)}}:\C{F}_{c_1(y)}\to \C{F}_{x}$. In particular,
for each $y\in c_2^{-1}(x)\cap c_1^{-1}(x)$, we get an endomorphism 
$u_y:\C{F}_x\to\C{F}_x$.
\end{Emp}

\begin{Def}
Let $c:C\to X\times X$ be a correspondence. We say that a closed subset $Z\subset X$ is 
{\em locally $c$-invariant}, if for each $x\in Z$ there exists an open neighborhood 
$U\subset X$ of 
$x$ such that $c_1(c_2^{-1}(Z\cap U))$ is set-theoretically 
contained in $Z\cup (X\sm U)$. 
\end{Def}


\subsection{Main theorem}

\begin*
\vskip 8truept
\end*

\begin{Not} \label{N:degree}
Let $f:Y\to X$ be a morphism of Noetherian schemes. 

a) For a closed reduced subscheme $Z\subset X$, we denote by $ram(f,Z)$ 
({\em the ramification of $f$ at $Z$}) the smallest positive integer $m$
such that $(\C{I}_{(f^{-1}(Z))_{red}})^m\subset\C{I}_{f^{-1}(Z)}$.

b) If $f$ is quasi-finite, we denote by $ram(f)$  ({\em the ramification degree of $f$}) 
the maximum of $ram(f,x)$, where $x$ runs over the set of all closed points of $X$.
\end{Not}

Now we are ready to formulate our  main result.

\begin{Thm} \label{T:Del} 
Let $c:C\to X\times X$ be a correspondence defined over $\fq$.

a) Assume that $c_2$ is quasi-finite. Then for every $n\in\B{N}$ with $q^n>ram(c_2)$, 
the set $Fix(c^{(n)})$ is finite. 

b) Let  $U\subset X$ be an open subset defined over $\fq$ such that 
$c_1|_{c_1^{-1}(U)}$ is proper, $c_2|_{c_2^{-1}(U)}$ is quasi-finite,
and $X\sm U$ is locally $c$-invariant.

Then there exists a positive integer $d\geq ram(c_2|_{c_2^{-1}(U)})$ with 
the following property: For every  
$\C{F}\in D_{ctf}^b(X,\La)$ with $\C{F}|_{X\sm U}=0$, 
every $n\in\B{N}$ with $q^n>d$ and every  $c^{(n)}$-morphism 
$u: c_{2!}(c^{(n)}_1)^*\C{F}\to\C{F}$, we have an equality
\begin{equation} \label{Eq:fuj}
\Tr(R\Gm_c(u))=\sum_{y\in Fix(c^{(n)})\cap c'^{-1}(U)}\Tr(u_y).
\end{equation}

c) In the notation of b), assume that $X$ and $C$ are proper. Then \\  
$d:=max\{ram(c_2|_{c_2^{-1}(U)}),ram(c_2,X\sm U)\}$ satisfies the conclusion
of b). 
\end{Thm}

\begin{Rem} 
a) Note that both sides of \form{fuj} are well-defined. Namely, 
$R\Gm_c(u)$ was defined in \re{rgm} b), $u_y$ was defined in \re{naive}, 
and the sum is finite by a).

b) In the notation of \ref{T:Del} b), assume that $\C{F}\in D_{ctf}^b(X,\La)$ is 
equipped with a morphism $\psi:\Fr_X^*\C{F}\to\C{F}$ (say, $\C{F}$ is a Weil sheaf) 
and with a $c$-morphism $u:c_{2!}c_1^*\C{F}\to\C{F}$.  Then for each $n\in\B{N}$,
 $\C{F}$ is equipped with a $c^{(n)}$-morphism 
$u^{(n)}:=u\circ \psi^n:c_{2!}(c^{(n)}_1)^*\C{F}\to\C{F}$, so one can apply formula 
\form{fuj}. In the case $U=X$, the assertion thus reduces to Deligne's conjecture proven by 
Fujiwara \cite{Fu}.

c) The constant $d$ in \ref{T:Del} b) can be made explicit. 
Namely, one can see from the proof that the picture can be compactified, and then 
\ref{T:Del} c) gives an estimate for $d$.
\end{Rem}

\section{Proof of the Main Theorem}

\subsection{Push-forward of cohomological correspondences}

\begin*
\vskip 8truept
\end*

\begin{Def} \label{D:mor} 
By a {\em morphism} from a correspondence  $c:C\to X_1\times X_2$ to a correspondence 
$b:B\to Y_1\times Y_2$ we mean a triple $[f]=(f_1,f^{\diez},f_2)$ making the following diagram commutative 
\begin{equation} \label{Eq:funct}
\CD 
        X_1   @<{c_1}<<                    C        @>{c_2}>>         X_2\\
        @V{f_1}VV                        @V{f^{\diez}}VV                       @VV{f_2}V\\
        Y_1 @<{b_1}<<                    B    @>{b_2}>>            Y_2.
\endCD 
\end{equation}
\end{Def}

\begin{Emp} \label{E:pushf}
{\bf Construction.}
In the notation of \ref{D:mor} assume that either 

(i) the left hand square of \form{funct} is Cartesian, or 

(ii) morphisms $f_1$ and $f^{\diez}$ are proper, or 

(iii) morphisms $c_1$ and $b_1$ are proper.

In all these cases, we have base change morphisms 
$BC:b_1^* f_{1!}\to f^{\diez}_! c_1^*$.
Hence every $c$-morphism $u:c_{2!}c_1^*\C{F}_1\to\C{F}_2$
gives rise to a $b$-morphism
\[
[f]_!(u):b_{2!}b_1^*(f_{1!}\C{F}_1)\overset{BC}{\lra} b_{2!}f^{\diez}_! c_1^*\C{F}_1=
 f_{2!}c_{2!}c_1^*\C{F}_1\overset{u}{\lra} f_{2!}\C{F}_2.
\]
\end{Emp}

\begin{Emp} \label{E:rgm'}
{\bf Example.} For each correspondence $c$, there is a structure morphism 
$[\pi]_c$ from $c$ to the trivial correspondence $c_{tr}:\Spec k\to\Spec k\times\Spec k$.
Moreover, if  $c_1$ is proper, then $[\pi]_c$ satisfies assumption (iii) of 
\re{pushf}, and $[\pi]_{c!}$ coincides with 
the map $R\Gm_c$, defined in \re{rgm} a).
\end{Emp}

\subsection{Trace maps and local terms} 

\begin*
\vskip 8truept
\end*
 
\begin{Emp} \label{E:locterms}
{\bf Definition of trace maps.}
Let $c:C\to X\times X$ be a correspondence.

a) Denote by 
$\un{\C{Tr}}:\C{RHom}(c_1^*\C{F},c_2^!\C{F})\to\Dt'_*K_{Fix(c)}$
the composition
\[
\C{RHom}(c_1^*\C{F},c_2^!\C{F})\isom
c^!(\B{D}\C{F}\pp\C{F})\to c^!\Dt_*K_X\isom\Dt'_* c'^!K_X=\Dt'_*K_{Fix(c)},
\]
where the first map is defined in \cite[(3.1.1) and (3.2.1)]{Il} using the standard 
identification $\C{RHom}(\C{A},\C{B})\cong\Dt^!(\B{D}\C{A}\pp\C{B})$, the second one is 
induced by the map $\B{D}\C{F}\pp\C{F}\to\Dt_*K_{X}$, adjoint to the evaluation map
$\Dt^*(\B{D}\C{F}\pp\C{F})=\B{D}\C{F}\otimes\C{F}\to K_{X}$,
and the last one is the base change isomorphism  
$c^!\Dt_*\isom \Dt'_* c'^!$.

Since  
$H^0(C,\C{RHom}(c_1^*\C{F},c_2^!\C{F}))=\Hom(c_1^*\C{F},c_2^!\C{F})\cong \Hom(c_{2!}c_1^*\C{F},\C{F})$,
the map $H^0(C,\un{\C{Tr}})$ gives rise to a map
\begin{equation} \label{Eq:lterm}
\C{Tr}=\C{Tr}_c:\Hom(c_{2!}c_1^*\C{F},\C{F})\to H^0(C,\Dt'_*K_{Fix(c)})=H^0(Fix(c),K_{Fix(c)}),
\end{equation}
which we call the {\em trace map}.

b) For an open subset $\beta$ of $Fix(c)$, we denote by 
\[
\C{Tr}_{\beta}:\Hom(c_{2!}c_1^*\C{F},\C{F})\to H^0(\beta, K_{\beta})
\]
the composition of $\C{Tr}$ and the restriction map $H^0(Fix(c),K_{Fix(c)})\to  H^0(\beta, K_{\beta})$.

 If, moreover, $\beta$ is proper over $k$, we denote by 
$LT_{\beta}:Hom(c_{2!}c_1^*\C{F},\C{F})\to\La$
the composition of $\C{Tr}_{\beta}$ and the integration map 
$H^0(\beta, K_{\beta})\to \La$.

If $\beta$ is a connected component of $Fix(c)$ which is proper over $k$, 
then $LT_{\beta}(u)$ is usually called the {\em local term} of $u$ at $\beta$.
\end{Emp}

\begin{Emp} \label{E:triv}
{\bf Example.} If $c=c_{tr}$ (see \re{rgm'}), then $H^0(Fix(c),K_{Fix(c)})=\La$, each 
$\C{F}$ is just a bounded complex of finitely generated free $\La$-modules (modulo homotopy), 
and the trace map $\C{Tr}_{c_{tr}}$ coincides with the usual trace map 
$\Hom(\C{F},\C{F})\to\La$.
\end{Emp} 

\begin{Rem}
Our trace map is equivalent to the map $\lan\cdot , \Id\ran$,
where 
\[
\lan\cdot,\cdot\ran:\Hom(c_{2!}c_1^*\C{F},\C{F})\otimes\Hom(\C{F},\C{F})\to 
H^0(\beta, K_{\beta})
\]
is the pairing, associated by Illusie \cite[(4.2.5)]{Il}, to a pair of correspondences 
$c:C\to X\times X$ and $\Dt:X\to X\times X$. However, our notion is more 
elementary. 
\end{Rem}

As in \cite[Cor. 4.5]{Il}, the trace maps commute with proper push-forwards. 

\begin{Prop} \label{P:pushf} 
Let $[f]=(f,f^{\diez},f)$ be a morphism from a correspondence $c:C\to X\times X$ to 
$b:B\to Y\times Y$ such that $f$ and $f^{\diez}$ are proper. 
Then the morphism $f':Fix(c)\to Fix(b)$ induced by $[f]$ is proper as well, 
and for every $c$-morphism
$u:c_{2!}c_1^*\C{F}\to\C{F}$, we have an equality
\begin{equation} \label{Eq:pushf}
\C{Tr}_b([f]_!(u))=f'_!(\C{Tr}_c(u))\in H^0(Fix(b),K_{Fix(b)}).
\end{equation}
\end{Prop}

Applying the proposition to the case when $X$ and $C$ are proper over $k$ and 
$[f]$ is the structure morphism $[\pi]_c$ of \re{rgm'}, we deduce the 
well-known Lefschetz--Verdier trace formula (\cite[Cor. 4.7]{Il}).

\begin{Cor} \label{C:LTF}
Let  $c:C\to X\times X$ be a correspondence 
such that $C$ and $X$ are proper over $k$. Then for every $c$-morphism 
 $u: c_{2!}c_1^*\C{F}\to\C{F}$, we have an equality
\begin{equation} \label{Eq:ltf} 
\Tr(R\Gm_c(u))=\sum_{\beta\in \pi_0(Fix(c))}LT_{\beta}(u).
\end{equation}
\end{Cor}

\subsection{Invariant subsets}

\begin*
\vskip 8truept
\end*

\begin{Def}
Let $c:C\to X\times X$ be a correspondence, and $Z\subset X$ a closed subset.

a) We say that $Z$ is {\em $c$-invariant}, if $c_1(c_2^{-1}(Z))$ is set-theoretically contained in $Z$. 

b) We say that $Z$ is {\em $c$-invariant in a neighborhood of fixed points},
if there exists an open neighborhood $W\subset C$ of $Fix(c)$ such that
$Z$ is $c|_W$-invariant, where $c|_W:W\to X\times X$ is the restriction of $c$.
\end{Def}

\begin{Lem} \label{L:invneib}
Let $c:C\to X\times X$ be a correspondence, and $Z\subset X$ a closed subset.

a) The set $W(Z):=C\sm [\ov{c_2^{-1}(Z)\sm c_1^{-1}(Z)}]$ is the largest open subset 
$W\subset C$ such that $Z$ is $c|_W$-invariant, where by $\ov{\cdot}$ we denote the closure.

b) $Z$ is locally $c$-invariant if and only if for each irreducible component 
$S$ of  $c_2^{-1}(Z)\sm c_1^{-1}(Z)$, the closures of $c_1(S)$ and  $c_2(S)$ in $X$
do not intersect.

c) If $Z$ is locally $c$-invariant, then $Z$ is  $c$-invariant in a neighborhood of fixed points.
\end{Lem}

\begin{proof}
a) and b) follow from definitions, c) follows from a) and b).
\end{proof}

\begin{Emp} \label{E:qfinite} 
{\bf Example.} If $c_2$ is quasi-finite, then every  closed point $x\in X$ is 
locally $c$-invariant. Indeed, $U:=X\sm [c_1(c_2^{-1}(x))\sm x]$ is the required 
open neighborhood. As a result, every closed point $x\in X$ is 
locally $c$-invariant in a neighborhood of fixed points.
\end{Emp}

\begin{Def} 
Let $c:C\to X\times X$ be a correspondence. We say that a correspondence 
$\ov{c}:\ov{C}\to \ov{X}\times\ov{X}$ is a {\em compactification} of $c$, 
if $\ov{X}$ and $\ov{C}$ are proper over $k$, $C\subset\ov{C}$ and 
$X\subset\ov{X}$ are open subsets, and $c$ is the restriction of $\ov{c}$.
\end{Def}

The following lemma can be deduced from \rl{invneib} b).

\begin{Lem} \label{L:comp}
Let $c:C\to X\times X$ be a correspondence and $U\subset X$ an open subset such that
$c_1^{-1}(U)$ is dense in $C$, $c_1|_{c_1^{-1}(U)}$ is proper and  
$X\sm U$ is locally $c$-invariant.
Then there exists a compactification  $\ov{c}:\ov{C}\to \ov{X}\times \ov{X}$ of 
$c$ such that $\ov{X}\sm U$ is locally $\ov{c}$-invariant.
\end{Lem}

\begin{Emp} \label{E:restr2}
{\bf Restriction of correspondences.}
Let $c:C\to X\times X$ be a correspondence,
 $u:c_{2!}c_1^*\C{F}\to\C{F}$ a $c$-morphism, and
$Z\subset X$ a closed subset.

a) If $Z$ is $c$-invariant, then $c$ induces a correspondence
$c|_Z:c_2^{-1}(Z)_{red}\to Z\times Z$, hence by \re{restr}, $u$ restricts to a 
$c|_Z$-morphism $u|_Z$. 

b) In general, let $W=W(Z)\subset C$ be as in \rl{invneib} a). Then $Z$ is $c|_W$-invariant, 
and we denote the correspondence $(c|_W)|_Z$ defined in a) simply by $c|_Z$. 
Moreover, by \re{restr}, $u$ restricts to a $c|_W$-morphism $u|_W$, 
hence by a) to a  $c|_Z$-morphism $u|_Z:=(u|_W)|_Z$.
\end{Emp}

\begin{Emp} \label{E:naive2}
{\bf Example.}
If $c_2$ is quasi-finite, and $Z$ is a closed point $x$, then 
$c|_Z=c|_x$ is the correspondence $c_1^{-1}(x)\cap c_2^{-1}(x)\to\{x\}\times\{x\}$.
Moreover, for each $y\in c_1^{-1}(x)\cap c_2^{-1}(x)$, the restriction of $u|_x$ to 
$\{y\}\to\{x\}\times \{x\}$ equals the endomorphism $u_y:\C{F}_x\to\C{F}_x$ defined in 
\re{naive}. Using \re{triv} we see that $LT_y(u|_x)=\Tr(u_y)$.
\end{Emp}

\subsection{Contracting correspondences}

\begin*
\vskip 8truept
\end*

\begin{Def} \label{D:contr}
Let $c:C\to X\times X$ be a correspondence, and $Z\subset X$ a closed subscheme.

a) We say that $c$ {\em stabilizes} $Z$,  if 
$c_1(c_2^{-1}(Z))$ is scheme-theoretically contained in $Z$, i.e.,
if $c_1^{\cdot}(\C{I}_{Z})\subset c_2^{\cdot}(\C{I}_{Z})\cdot\C{O}_C$.

b) We say that $c$ is {\em contracting near} $Z$, if $c$ stabilizers $Z$ 
and there exists $n\in \B{N}$ such that 
$c_1^{\cdot}(\C{I}_{Z})^n\subset c_2^{\cdot}(\C{I}_{Z})^{n+1}\cdot\C{O}_{C}$.

c) We say that $c$ is {\em contracting near $Z$ in a neighborhood of fixed points},
if there exists an open neighborhood $W\subset C$ of $Fix(c)$ such that
$c|_W:W\to X\times X$ is contracting near $Z$.
\end{Def}

\begin{Rem} \label{R:contr}
a) A geometric characterization of a contracting correspondence 
will be given in \rr{contr2}.

b) If a correspondence $c$ is contracting near
$Z$, then $c^{rig}$ is contracting near $Z^{rig}$ in the sense of \cite[Def. 3.1.1]{Fu}.
Furthermore, it is likely that the two notions are equivalent.
\end{Rem} 
 
The proof of the following crucial result will occupy \rs{contrcor}.

\begin{Thm} \label{T:locterms}
Let $c:C\to X\times X$ be a correspondence contracting near a closed subscheme 
$Z\subset X$ in a neighborhood of fixed points, and let $\beta$ be an open connected subset of 
$Fix(c)$ such that  $c'(\beta)\cap Z\neq\emptyset$. Then 

a) $c'(\beta)$ is contained in $Z$, hence $\beta$ 
is an open connected subset of $Fix(c|_Z)$.

b) For every $c$-morphism $u:c_{2!}c_1^*\C{F}\to\C{F}$, we have
$\C{Tr}_{\beta}(u)=\C{Tr}_{\beta}(u|_Z)$. In particular, if $\beta$ is proper, then 
$LT_{\beta}(u)=LT_{\beta}(u|_Z)$.
\end{Thm}

To apply the result, we will use the following lemma.

\begin{Lem} \label{L:twisting}
Let $c:C\to X\times X$ be a correspondence defined over $\fq$.

a) Let $Z$ be a locally $c$-invariant closed subset of $X$ defined over $\fq$.   
Then for each  $n\in\B{N}$ with $q^n>ram(c_2,Z)$, the correspondence $c^{(n)}$ is contracting 
near $Z$ in a neighborhood of fixed points.

b) If $c_2$ is quasi-finite, then for each $q^n>ram(c_2)$, the correspondence $c^{(n)}$ 
is contracting near every closed point $x$ of $X$ 
in a neighborhood of fixed points.
\end{Lem}

\subsection{Proof of \rt{Del}.}

\begin*
\vskip 8truept
\end*

\begin{Emp}
{\bf Proof of a).} 
The assertion follows from \rl{twisting} b), \rt{locterms} a) and the fact that
for every closed point $x$ of $X$, the set $Fix(c|_x)$ is finite.
\end{Emp}

\begin{Emp}
{\bf Proof of c).} 
By \rco{LTF}, we have an equality
\begin{equation} \label{Eq:lef}
\Tr(R\Gm_c(u))=\sum_{\beta\in\pi_0(Fix(c^{(n)}))}LT_{\beta}(u).
\end{equation}
Pick any $\beta\in \pi_0(Fix(c^{(n)}))$ such that $c_2(\beta)\subset X\sm U$. 
Since $q^n>ram(c_2,X\sm U)$, we conclude from \rl{twisting} a) and \rt{locterms}, that 
$\beta$ is a connected component of $Fix(c|_{X\sm U}^{(n)})$, and  
$LT_{\beta}(u)$ equals $LT_{\beta}(u|_{X\sm U})$. 
As $\C{F}|_{X\sm U}=0$, we get that  $LT_{\beta}(u|_{X\sm U})=0$, 
hence also $LT_{\beta}(u)=0$.

Pick now any  $\beta\in \pi_0(Fix(c^{(n)}))$ such that 
$c_2(\beta)\cap U\neq\emptyset$. As  $q^n>ram(c_2|_{c_2^{-1}(U)})$, 
such $\beta$ is simply a point $y\in Fix(c^{(n)})\cap c'^{-1}(U)$ (by a)). 
Moreover, by \rl{twisting} b), \rt{locterms} b) 
and \re{naive2}, $LT_{\beta}(u)$ equals  $LT_y(u|_x)=Tr(u_y)$.
This shows that the right hand side of \form{lef} is equal to that of \form{fuj}, as claimed.
\end{Emp}

\begin{Emp}
{\bf Proof of b).} 
For the proof we can replace $c$ and $u$ by their restrictions to  $c_1^{-1}(U)$. 
Then assumptions of \rl{comp} are satisfied, hence there exists a compactification 
$\ov{c}:\ov{C}\to \ov{X}\times \ov{X}$ of $c$ such that $\ov{X}\sm U$ is 
locally $\ov{c}$-invariant. We claim that 
$d:=max\{ram(c_2|_{c_2^{-1}(U)}),ram(\ov{c}_2,\ov{X}\sm U)$  
has the required property.

Let $c^0:C^0=c_1^{-1}(U)\to U\times X$ be the restrictions of $c$, $u^0$ 
the restrictions of $u$ to $c^0$, and $[j]=(j^0,j_{C^0},j)$ the inclusion map of 
$c^0$ into $\ov{c}$. Then $[j]$ satisfies assumption (iii) of \re{pushf}, therefore $u^0$ 
extends to a $\ov{c}$-morphism 
$\ov{u}:=[j]_!(u^0):\ov{c}_{2!}\ov{c}_1^*\ov{\C{F}}\to\ov{\C{F}}$, where 
$\ov{\C{F}}:=j_!(\C{F})=j^0_!(\C{F}|_U)$. 

We claim that the equality \form{fuj} for $c, U$ and $u$ is equivalent to that for 
$\ov{c}, U$ and $\ov{u}$. Indeed, since $[\pi]_{c^0}=[\pi]_{\ov{c}}\circ [j]$, 
the equality $\Tr(R\Gm_c(\ov{u}))=\Tr(R\Gm_c(u))$ follows from \re{rgm'}, 
while the equality $\Tr(\ov{u}_y)=\Tr(u_y)$ for all $y\in Fix(c^{(n)})\cap c'^{-1}(U)$
is clear. Now the assertion follows from part c) of the theorem proven above.
\end{Emp}

\section{Local terms for contracting correspondences.} \label{S:contrcor}

This section is devoted to the proof of \rt{locterms}. 

\subsection{Additivity of the trace maps} 

\begin*
\vskip 8truept
\end*

\begin{Not} \label{N:clop}
Let $c:C\to X\times X$ be a correspondence, $u:c_{2!}c_1^*\C{F}\to \C{F}$ a 
$c$-morphism, and  $Z\subset X$ a $c$-invariant closed subset.

a) Let $c|_Z$ and $u|_Z$ be as in \re{restr2} a), and let
$[i]_Z$ be the closed embedding of $c|_Z$ into $c$. Then $[i]_Z$ satisfies assumption  
(ii) of \re{pushf}, hence $u$ gives rise to a $c$-morphism 
$[i]_{Z!}(u|_Z)$.

b) Put $U:=X\sm Z$. Then $c_1^{-1}(U)\subset c_2^{-1}(U)$, 
and we denote by  $c|_U:c_1^{-1}(U)\to U\times U$ the restriction of $c$, by  $u|_U$
the restriction of $u$ to $c|_U$ (see \ref{E:restr}), and by  
$[j]_U$ the open embedding of $c|_U$ into $c$. Since $[j]_U$ satisfies assumption (i) of   
\re{pushf}, $u|_U$ extends to a $c$-morphism  
$[j]_{U!}(u|_U)$.
\end{Not}

As in \cite[Prop. 2.4.3]{Pi} and \cite[4.13]{Il}, the trace maps are additive.

\begin{Prop} \label{P:add}
In the notation of \ref{N:clop}, we have an equality
\[
\C{Tr}_c(u)= \C{Tr}_c([i]_{Z!}(u|_Z))+\C{Tr}_c([j]_{U!}(u|_U)).
\]  
\end{Prop}

\subsection{Specialization} \label{S:spec}

\begin*
\vskip 8truept
\end*

\begin{Not}
For a scheme $X$ over $k$, set $X_{\B{A}}:=X\times \B{A}^1$.
For a morphism $f:X\to Y$ of schemes over $k$, set
$f_{\B{A}}:=f\times\Id_{\B{A}^1}:X_{\B{A}}\to Y_{\B{A}}$. For a  
scheme $\wt{X}$ over $\B{A}^1$, we denote by $\wt{X}_s$ its fiber over $0\in\B{A}^1$, 
and denote by
$\Psi_{\wt{X}}: D_{ctf}^b(\wt{X},\La)\to D_{ctf}^b(\wt{X}_{s},\La)$ the corresponding 
functor of nearby cycles. 
\end{Not}

\begin{Emp} \label{E:specfunc}
{\bf Specialization functor.}

a) We say that a scheme $\wt{X}$ over $\B{A}^1$ {\em lifts} a scheme $X$ over $k$, 
if it is equipped with a morphism $\varphi=\varphi_X:\wt{X}\to X$ such 
the corresponding morphism $\wt{X}\to X_{\B{A}}=X\times\B{A}^1$ is an isomorphism 
over $\B{A}^1\sm\{0\}$. In this case, we define a functor 
\[
sp_{\wt{X}}:=\Psi_{\wt{X}}\circ \varphi^*:D_{ctf}^b(X,\La)\to D_{ctf}^b(\wt{X}_{s},\La).
\]

b)  We say that a morphism $\wt{f}:\wt{X}\to\wt{Y}$ of schemes over $\B{A}^1$
{\em lifts} a  morphism $f:X\to Y$ of schemes over $k$, if $\wt{X}$ lifts $X$, $\wt{Y}$ lifts $Y$,
and $\varphi_Y\circ\wt{f}=f\circ\varphi_X$. 
In this case, we have base change morphisms
$BC^*:\wt{f}_s^* sp_{\wt{Y}}\to sp_{\wt{X}} f^*$, 
$BC_*: sp_{\wt{Y}} f_*\to \wt{f}_{s*} sp_{\wt{X}}$, 
$BC^!: sp_{\wt{X}} f^!\to \wt{f}_s^! sp_{\wt{Y}}$ and
$BC_!:\wt{f}_{s!} sp_{\wt{X}}\to sp_{\wt{Y}} f_!$, induced by the corresponding
base change morphisms for $\Psi$.
\end{Emp}

\begin{Emp} \label{E:trivdef}
{\bf Examples.}
a) If $\wt{X}=X_{\B{A}}$ and $\varphi$ is the projection map, then $\wt{X}_s=X$, and the functor 
$sp_{\wt{X}}$ is isomorphic to the identity functor. 

b) If $f:X\to\Spec k$ and $\wt{f}:\wt{X}\to\B{A}^1$ are the structure morphisms, 
then the composition $BC^!\circ BC_*:f_* f^!\La=sp_{\B{A}^1}f_* f^!\La\to 
\wt{f}_{s*}\wt{f}_s^! sp_{\B{A}^1}\La=\wt{f}_{s*}\wt{f}_s^!\La$ defines a morphism 
$H^0(X,K_X)\to H^0(\wt{X}_s,K_{\wt{X}_s})$, which we will denote simply by 
$sp_{\wt{X}}$.
\end{Emp}

\begin{Emp} \label{E:cor}
{\bf Specialization of correspondences.}
Let  $c:C\to X\times X$ be a correspondence, $u:c_{2!}c_1^*\C{F}\to\C{F}$
a $c$-morphism, 
$\wt{c}:\wt{C}\to \wt{X}\times \wt{X}$ a correspondence over $\B{A}^1$ lifting $c$, and 
$\wt{c}_{s}$ the fiber of $\wt{c}$ over $0\in\B{A}^1$. Then $u$ 
gives rise to a $\wt{c}_{s}$-morphism 
\[
sp_{\wt{c}}(u):\wt{c}_{s2!}\wt{c}_{s1}^* sp_{\wt{X}}\C{F}
\overset{BC^*}{\lra}\wt{c}_{s2!} sp_{\wt{C}} c_{1}^*\C{F}
\overset{BC_!}{\lra} sp_{\wt{X}}c_{2!}c_1^*\C{F}
\overset{u}{\lra}
 sp_{\wt{X}}\C{F}.
\] 
\end{Emp}

As in \cite[Prop. 1.7.1]{Fu}, the trace maps commute with the specialization. 

\begin{Prop} \label{P:spec}
In the notation of \re{cor}, assume that $\wt{c}$ extends to a lifting of some 
compactification  $\ov{c}:\ov{C}\to \ov{X}\times \ov{X}$ of $c$. Then 
we have an equality
\begin{equation} \label{Eq:spec}
\C{Tr}_{\wt{c}_s}(sp_{\wt{c}}(u))=sp_{Fix(\wt{c})}(\C{Tr}_c(u))
\in H^0(Fix(\wt{c}_s),K_{Fix(\wt{c}_s)}).
\end{equation}
\end{Prop}


\subsection{Deformation to the normal cone} \label{SS:normcone}

\begin*
\vskip 8truept
\end*

We will apply the specialization in the following particular case.
 
\begin{Not} \label{N:deform}
Let $X$ be a scheme over $k$, $Z\subset X$ a closed subscheme.

a) Denote by  $\wt{X}_{Z}$ the spectrum of the 
$\C{O}_{X}$-subalgebra
$\C{O}_{X}[t,\frac{\C{I}_Z}{t}]\subset\C{O}_{X}[t,t^{-1}]$. The embedding
$\C{O}_{X}[t]\hra\C{O}_{X}[t,\frac{\C{I}_Z}{t}]$ gives rise to the birational projection
$\wt{\varphi}:\wt{X}_{Z}\to X_{\B{A}}$, which is an isomorphism over $\B{A}^1\sm\{0\}$. 

b) The special fiber $(\wt{X}_{Z})_s=\C{Spec}(\bigoplus_{n=0}^{\infty}
(\C{I}_Z)^n/(\C{I}_Z)^{n+1})$ is the normal cone of $X$ to $Z$, which we denote by 
$N_{Z}(X)$. 

c) The projection $\C{O}_{X}[t,\frac{\C{I}_Z}{t}]\to
(\C{O}_{X}/\C{I}_{Z})[t]=\C{O}_{Z}[t]$ defines a closed embedding 
$\wt{\iota}:Z_{\B{A}}\hra\wt{X}_{Z}$. The special fiber
$\wt{\iota}_s:Z\hra N_{Z}(X)$  identifies $Z$ with the zero section of  $N_{Z}(X)$.

d) Since $sp_{Z_{\B{A}}}$ is the identity functor (see \re{trivdef}),  
the map $BC^*$ (from \re{specfunc}) for the embedding $\wt{\iota}$ from c) defines a 
morphism
\begin{equation} \label{Eq:Ve} 
sp_{\wt{X}_{Z}}(\C{F})|_Z\to\C{F}|_Z. 
\end{equation}
\end{Not}

The following property of the deformation to the normal cone, 
proven in \cite[$\S$8, (SP5)]{Ve}, is crucial for the whole proof.

\begin{Lem} \label{L:verdier}
The morphism \form{Ve} is an isomorphism.
\end{Lem}

The following lemma is straightforward.

\begin{Lem} \label{L:normfunc}
Let  $f:X_1\to X_2$ be a morphism of schemes over $k$, $Z_2\subset X_2$ a closed subscheme,
and $Z_1$ a closed subscheme of $f^{-1}(Z_2)$, that is, 
$f^{\cdot}(\C{I}_{Z_2})\subset\C{I}_{Z_1}$. 

a)  The morphism $f_{\B{A}}:X_{1\B{A}}\to X_{2\B{A}}$ lifts uniquely  
to the map $\wt{f}:(\wt{X_1})_{Z_1}\to(\wt{X_2})_{Z_2}$.

b) The image $\wt{f}_s({N}_{Z_1}(X_1))$ is supported set-theoretically at  
the zero section of $N_{Z_2}(X_2)$ if and only if there exists $n\in \B{N}$ such that 
$f^{\cdot}(\C{I}_{Z_2})^n\subset (\C{I}_{Z_1})^{n+1}$.
\end{Lem}

\begin{Emp} \label{E:defcor}
{\bf Deformation of correspondences.}
Let $c:C\to X\times X$ be a correspondence, and  $Z\subset X$ a closed subscheme. Then by 
\rl{normfunc}, $c$ lifts to a correspondence 
$\wt{c}_{Z}:\wt{C}_{c^{-1}(Z\times Z)}\to\wt{X}_{Z}\times\wt{X}_{Z}$. 
Moreover, if $\ov{c}:\ov{C}\to\ov{X}\times\ov{X}$ is a compactification of $c$, and 
$\ov{Z}\subset\ov{X}$ is the closure of $Z$, then correspondence 
$\wt{\ov{c}}_{\ov{Z}}$ extends $\wt{c}_{Z}$ and lifts $\ov{c}$. 
In particular, \rp{spec} holds for $c$ and $\wt{c}_Z$.
\end{Emp}

\begin{Rem} \label{R:contr2}
Recall that a correspondence $c:C\to X\times X$ stabilizes a closed subscheme $Z\subset X$ 
if and only if $c^{-1}(Z\times Z)=c_2^{-1}(Z)$. Therefore by \rl{normfunc} b), $c$ is 
contracting near $Z$ if and only if $c$  stabilizes $Z$ and the image of 
$(\wt{c}_Z)_{1s}$ is 
supported  set-theoretically at the zero section $Z\subset N_{Z}(X)$. 
\end{Rem}

\subsection{Proof of \rt{locterms}}

\begin*
\vskip 8truept
\end*

Choose  an open neighborhood $W\subset C$ of $Fix(c)$ such that $c|_W$ is contracting near $Z$.
Replacing $c$ by $c|_W$, we can assume that $c$ is contracting near $Z$.
Moreover, replacing further $C$ by an open subset $C\sm [Fix(c)\sm \beta]$ we can assume that
$Fix(c)=\beta$, hence $\C{Tr}_{\beta}=\C{Tr}_c$. For the proof we apply the 
construction of \re{defcor}.

\begin{Emp}
{\bf Proof of a).}
By \rr{contr2}, the image of $(\wt{c}_{Z})_{1s}$ is supported at 
$Z\subset N_{Z}(X)$. Since $\wt{Fix(c)}_{c'^{-1}(Z)}$ is a closed subscheme of 
$Fix(\wt{c}_{Z})$, the image of $\wt{c}'_s:N_{c'^{-1}(Z)}(\beta)\to  N_{Z}(X)$ is then supported 
at $Z$ as well.
Using the ``only if'' statement of \rl{normfunc} b) we conclude that
$\beta_{red}=(c')^{-1}(Z)_{red}$, as claimed.
\end{Emp}

\begin{Emp}
{\bf Proof of b).}
Put $U:=X\sm Z$. Then we have 
$\C{Tr}_c(u)=\C{Tr}_c([i]_{Z!}(u|_Z))+\C{Tr}_c([j]_{U!}(u|_U))$ (by \rp{add}) and   
$\C{Tr}_c([i]_{Z!}(u|_Z))=i'_{Z!}\C{Tr}_{c|_Z}(u|_Z)$ (by \rp{pushf}). Since $i'_{Z!}$ is the 
identity map (by a)), it remains to prove that\\ $\C{Tr}_c([j]_{U!}(u|_U))=0$. 
For this, we may replace $\C{F}$ by $j_!(\C{F}|_U)$ and 
$u$ by $[j]_{U!}(u|_U)$. In this case $\C{F}|_Z=0$, and we claim that $\C{Tr}_c(u)=0$.

By \rp{spec} for correspondences $c$ and $\wt{c}_Z$, 
is will suffice to show that the specialization map $sp_{\wt{c}_Z}$
vanishes and the map $sp_{Fix(\wt{c}_Z)}$ is an isomorphism. 
Since the image of $(\wt{c}_Z)_{1s}$ is supported at $Z\subset N_{Z}(X)$ and 
$sp_{\wt{X}_Z}(\C{F})|_Z\cong \C{F}|_Z=0$ (by \rl{verdier}),
the sheaf $(\wt{c}_Z)_{1s}^* sp_{\wt{X}_Z}(\C{F})$ vanishes. This implying the
vanishing of $sp_{\wt{c}_{Z}}$. On the other hand, since $(\wt{c}_Z)_{2s}^{-1}(Z)$ equals 
$c_2^{-1}(Z)\subset N_{c_2^{-1}(Z)}(C)$, 
we see that $Fix(\wt{c}_Z)_{red}=(\beta_{\C{D}})_{red}\subset c_2^{-1}(Z)_{\C{D}}$. 
Thus  $sp_{Fix(\wt{c}_Z)}$ 
is an isomorphism by \re{trivdef}. 
\end{Emp} 

This completes the proof of \rt{locterms} and hence also of \rt{Del}.


\begin{thebibliography}{99}

\bibitem[Fu]{Fu}
K. Fujiwara, {\em Rigid geometry, Lefschetz--Verdier trace formula and Deligne's conjecture}, 
Invent. Math. {\bf 127} (1997), no. 3, 489--533.

\bibitem[KV]{KV}
D. Kazhdan and Y. Varshavsky, {\em On the cohomology of the moduli spaces of $F$-bundles: 
stable cuspidal Deligne--Lusztig part}, 
in preparation.

\bibitem[Il]{Il}
L. Illusie, {\em Formule de Lefschetz}, in {\em Coholologie $\ell$-adique et fonctions $L$}, SGA5, 
Lecture Notes in Mathematics {\bf 589}, Springer-Verlag, 1977, pp. 73--137.

\bibitem[Pi]{Pi}
R. Pink, {\em On the calculation of local terms in the Lefschetz--Verdier trace formula and 
its application to a conjecture of Deligne},  
Ann. of Math. {\bf 135} (1992), no. 3, 483--525.

\bibitem[Ve]{Ve}
J.-L. Verdier, {\em Sp\'ecialisation de faisceaux et monodromie mod\'er\'ee}, in 
{\em Analysis and topology on singular spaces}, 332--364, Astérisque {\bf 101-102}, 1983,
pp. 332--364.





































\end{thebibliography}
\end{document}